\documentclass[12pt, psamsfonts]{amsart}
\usepackage{amsmath}
\usepackage{amsthm}
\usepackage{amssymb}
\usepackage{amscd}
\usepackage{amsfonts}
\usepackage{amsbsy}
\usepackage{epsfig}
\usepackage{ulem}

\usepackage[colorlinks=true,linkcolor=blue, urlcolor=red, citecolor=blue, backref]{hyperref}	
\usepackage[latin1]{inputenc}      
\usepackage{tikz}
\newtheorem{theorem}{Theorem}

\newtheorem{remark}[theorem]{Remark}

\newtheorem{example}[theorem]{Example}
\def\e{\varepsilon}
\def\N{\mathbb N }

\begin{document}

\title[LEO and expansivity imply the periodic specification property]{The role of continuity and expansiveness on leo and periodic specification properties}

\author{Serge Troubetzkoy}
\address{Aix Marseille Univ, CNRS, Centrale Marseille, I2M, Marseille, France\linebreak
postal address: I2M, Luminy, Case 907, F-13288 Marseille Cedex 9, France}
\email{serge.troubetzkoy@univ-amu.fr}
\urladdr{www.i2m.univ-amu.fr/perso/serge.troubetzkoy/}

\author{Paulo Varandas}

\address{Paulo Varandas, CMUP and Departamento de Matem\'atica, Universidade Federal da Bahia\\
Av. Ademar de Barros s/n, 40170-110 Salvador, Brazil}
\email{paulo.varandas@ufba.br}
\urladdr{https://sites.google.com/view/paulovarandas/}

 \date{}

\thanks{ 
The authors are deeply grateful to A. Blokh and F. Przytycki for very useful comments on an early version of the paper.
PV  was partially supported by CMUP (UIDB/00144/2020), which is funded by FCT (Portugal) with national (MEC) and European structural funds through the programs FEDER, under the partnership agreement PT2020, and by 
 Funda\c c\~ao para a Ci\^encia e Tecnologia (FCT) - Portugal, through the grant CEECIND/03721/2017 of the Stimulus of Scientific Employment, Individual Support 2017 Call.
}

\subjclass[2000]{37E05, 37B40, 46B25, 46.3}
\keywords{expansiveness, specification, locally eventually onto, entropy}

\begin{abstract}
In this short note we prove that a continuous map  of a compact manifold 
which is locally eventually onto and is expansive satisfies
the periodic specification property. We also discuss the role of continuity as a key condition in the previous 
characterization.
We include several examples to illustrate the relation between these concepts. 
\end{abstract}

\maketitle

\section{Introduction}
There is a well known hierarchy of topological properties involving the topological indecomposability 
of a dynamical system, as transitivity, topological mixing, and the specification property, among many 
others. The relation between these and many others has been addressed by Akin, Auslander 
and Nagar \cite{AkAuNa}.
The aim of this short note is to complement the above results, and to highlight the relation between 
the locally eventually onto (a dynamical property stronger than topological mixing) 
and the specification properties, and to make explicit 
the role of continuity on such characterization.
The specification property was first introduced by Bowen \cite{Bo},
for
 a survey of specification-like properties we recommend the following  article \cite{KwLaOp}, while for a survey of mixing properties we recommend the article \cite{AkAuNa}.

First let us recall some well known results.
Blokh \cite{Blo} showed that for  a continuous map of the interval  $ [0,1]$ the periodic specification property is equivalent to topological mixing (see e.g., \cite{Bu,Ru}).  So, while for continuous interval maps the picture is very well understood and most concepts of topological chaoticity coincide, this is  no longer true for more general metric spaces or whenever continuity breaks down.
Yet, the situation is well understood in the case of one-dimensional branched manifolds, where there is a characterization of
transitive dynamics due to Blokh~\cite{Blo0}. In brief terms, he established the following spectral  decomposition theorem:
either a transitive map $f$ of a graph has periodic points and it can be decomposed into $n$ connected subgraphs 
with finite pairwise intersections which are cyclically permuted and $f^n$ has the specification property, 
or $f$ is aperiodic and it is just a cycle of $n$ circles with $f^n$ being 
an irrational rotation. 
We refer the reader to \cite{Blo0} for more details.

It is noticeable that while any locally eventually onto continuous map has dense periodic sets, 
it may not
have periodic points (cf.\ \cite[Theorem~2.30 and {Example}~2.31]{AkAuNa}). In particular, a  locally eventually onto continuous map need not satisfy the periodic
specification property. Two results complement this discussion. First, expansiveness play a key role 
to bridge between the specification and periodic specification properties: a topological dynamical 
system satisfying the specification property and whose natural extension is expansive satisfies the periodic specification property 
(see e.g., \cite[Lemma~6]{KwLaOp}). Second, Yan, Yin and Wang \cite[Theorem~3.1]{YaYiWa} constructed an example of a topological mixing subshift, hence expansive, which does not even 
have the specification property.

 The situation is well understood in the case of continuous, open and distance expanding maps on 
 compact metric spaces. Indeed, since any such map satisfies the shadowing property 
 and periodic points are dense in the non-wandering set, 
 these admit a spectral decomposition theorem (see \cite[Theorem~4.3.8]{PU}).
 Moreover, any such map is topologically mixing map if and only it is locally eventually onto.  
  We refer the reader to \cite[Sections~4.2 and 4.3]{PU} for more details.
Similar, but slightly weaker results are  known if we drop the openness assumption, instead assume shadowing \cite{KKO}.

In general, 
while  the locally eventually onto property  need not ensure the periodic specification property,
the following result shows that expansiveness can act as a sufficient condition for it. 
We refer the reader to Section~\ref{sec:def} for definitions.

 \begin{theorem}\label{thm0}
Let $X$ be a compact and connected topological manifold.
 If the topological dynamical system $(X,f)$ is locally eventually onto and expansive then it has the periodic specification property.
 \end{theorem}

This  result is no longer true if one replaces the condition of $X$ being a compact 
topological manifold by the assumption of being an arbitrary compact metric space. We refer the reader to 
Example~\ref{Feliks}, where we present an expansive and locally eventually onto map for which the 
periodic specification fails.

Note that the specification property is a topological invariant,  hence we can ask whether such a property holds for the continuous map $f$ on $(X,d)$ or on the metric space $(X,d')$, for a equivalent  metric $d'$.
In the case of  compact and connected topological manifolds, Coven and Reddy~\cite{CoRe} 
constructed adapted metrics, proving that every expansive dynamics is indeed expanding
with respect to some equivalent metric. In particular, Theorem~\ref{thm0} is a direct consequence 
of the previous discussion together with the following:

\medskip

 \begin{theorem}\label{thm1}
 Assume that the topological dynamical system $(X,f)$ is expanding and locally eventually onto. Then $(X,f)$ 
 has the periodic specification property.
 \end{theorem}

The latter suggests that the failure of periodic specification for distance expanding maps is essentially related to
the lack of periodic points (see Example~\ref{Feliks}),  the non-compactness of the phase space (see Example~\ref{non-compact}),
or that the dynamics is not mixing.
 
Given the previous result it is natural to ask whether any locally eventually onto continuous map satisfies the specification property. 

\begin{remark}
It is worth mentioning that the situation is clear for continuous interval maps. Indeed, combining \cite[Theorem~B]{CM} and Blokh's theorem (cf.\ \cite[Theorem~3.4]{Ru}), it follows that the locally eventually onto property implies on the following conditions, which, for interval maps, are equivalent: \begin{itemize}\item[(i)] $f^2$ is transitive, \item[(ii)] $f^n$ is transitive for every $n\ge 1$, \item[(iii)] $f$ is topologically mixing,\item[(iv)] $f$ satisfies the specification property.\end{itemize}While the converse holds in the case of piecewise monotone continuous interval maps (cf.\  \cite[Lemma~4.1]{CM}), it fails for general continuous interval maps. In particular there are continuous interval maps satisfying the specification property for which the locally eventually onto property fails  (see e.g.,  \cite[Example~3]{BM}).
\end{remark}

On the positive direction, we notice that the same strategy used in Blokh's theorem 
can be used for conformal-like maps. 

\begin{theorem}\label{prop:blokh}
Every locally eventually onto, continuous and 
conformal-like map  on a compact metric space 
satisfies the periodic specification property.
\end{theorem}

\begin{remark}  In the definition of topological dynamical system, the assumption that the metric space is complete cannot be removed.  Throughout $\mathbb N$ be the set of non-negative integers (hence containing 0).
There exists a metric space $X\subset\{0,1,2\}^{\mathbb N}$ such that the shift map
$(X,\sigma)$ is locally eventually onto, it is clearly expansive, but fails even to present periodic points \cite[{Example}~2.31]{AkAuNa}. 
\end{remark}

Our second goal concerns describing the consequences of discontinuities on 
locally eventually onto maps.  
This is a problem dual to the one considered by Buzzi ~\cite{Bu},   the study of the
specification property for piecewise monotone interval maps. 
In the case of piecewise monotone continuous interval maps $f$, the transitivity for $f^2$ ensures the following  
``almost" locally eventually onto property: for any open interval $A$ and any closed interval $J\subset (0,1)$ there exists 
$N\ge 1$ so that $f^N(A)\supset J$ (see \cite[Theorem~6]{BM}). However, while the key step in this argument explores the density of periodic points, the classical argument that ensures the density of periodic points for expanding maps
does not apply for transitive piecewise expanding interval maps given that dynamical balls may fail to grow to a 
large scale.

We shall focus on 
important classes of dynamical systems known as $\beta$-expansions and $\beta$-shifts
(see e.g.,  \cite{Bl}). These can be realized by geometric models in the interval;  
for each $\beta>1$, the $\beta$-map is the $C^\infty$-piecewise expanding interval map
$T_\beta :[0,1) \to [0,1)$ given by 
$$
T_\beta(x) = \beta x - \lfloor \beta x \rfloor.
$$  
However, while the previous map is always expansive, and Markov for a countable set of parameters, 
$T_\beta$ does not satisfy the specification property for Lebesgue almost every parameter $\beta>1$ (cf.\  \cite{Bu}).
A characterization of the set of the values of $\beta$ which lead to maps with specification can be found in \cite{Sc}.
The next result shows that continuity is essential in Theorem~\ref{thm1}.

\begin{theorem}\label{thm-beta}
For Lebesgue almost every $\beta \in (1,+\infty)$ the map $T_\beta$:
\begin{itemize}
\item[(i)] is locally eventually onto;
\item[(ii)] is expansive;
\item[(iii)] does not satisfy the specification property \cite{Bu}.
\end{itemize}
\end{theorem}

We complete this section with two final comments on the relation between the specification and the locally eventually onto properties
for continuous maps in more general metric spaces. 
While any  Anosov diffeomorphism 
satisfies the specification property (see e.g  \cite{KH}), every volume preserving Anosov diffeomorphism is clearly 
not locally eventually onto. Nevertheless, on the converse direction, locally eventually onto maps
displaying non-uniform expansion often satisfy some measure-theoretical forms of specification (we refer the reader to 
\cite{Ol,Va} for the precise formulations).   

\section{Definitions}\label{sec:def}
Let $(X,d)$ be a compact metric space, and $f: X \to X$ a continuous map. 
We refer to $(X,f)$ as a \textit{dynamical system}. 

The map $f$ is called
\textit{locally eventually onto} (\textit{LEO}) if for every nonempty open set $U$ there is an $n\in{\N} := \{0,1,2,\dots\}$ such that $f^n(U)=X$.

For  integers $a \ge b \ge 0$ let $f^{[a,b]}(x) := \{f^j(x): a \le j
\le b\}$.
 
 A family of orbit segments $\{f^{[a_j,b_j]}(x_j)\}_{j=1}^n$
is an \textit{$N$-spaced specification} if \\$a_i - b_{i-1} \ge N$ for $2 \le i \le n$.  

We say that a specification $\{f^{[a_j,b_j]} (x_j)\}_{j=1}^n$ is \textit{$\e$-shadowed} by $y  \in X$ if
$$d(f^k(y),f^k(x_i)) \le \e \text{ for } a_i \le k \le b_i \text{ and } 1 \le i \le n.$$
  
We say that $(X, f)$ has the \textit{specification property} if for any  $\e > 0$
there is a constant $N = N(\e)$  such that any $N$-spaced specification
$\{f^{[a_j,b_j]} (x_j)\}_{j=1}^n$ is $\e$-shadowed by some $y \in X$.
If additionally, $y$ can be chosen in such
a way that $f^{b_n-a_0 +N}(y) = y$ then $(X,f)$ has the \textit{periodic specification property}.

The dynamical system $(X,f)$   is \textit{positively expansive}  if there exists 
$\alpha >0$ , called \textit{expansivity constant} of $f$, such that if $x,y \in X$  and  
$x \ne y$, then for some $n \ge 0$, $d(f^nn(x),f^n(y))>\alpha$.

The dynamical system $(X,f)$   is \textit{expanding}  if there are 
constants $\lambda > 1$ and $\delta_0 > 0 $ such that, for all $x,y,z\in X$, 
\begin{enumerate}
\item $d(f(x),f(y)) \ge \lambda d(x,y)$ whenever $d(x,y) < \delta_0$ and
\item $B(x,\delta_0) \cap f^{-1}(z)$ is a singleton whenever $d(f(x),z) < \delta_0$.
\end{enumerate}
A dynamical system $(X,f)$ satisfying condition (1) if called a \textit{distance expanding map}.
In any compact metrizable space, a continuous transformation is expanding if and only if it is open,
i.e., maps open sets  to open sets,
and distance expanding (see \cite[Lemma~1]{CoRe}).  
In \cite{PU} the authors describe the dynamical properties of such maps and obtaining, in particular,
density of periodic points, the shadowing property and a spectral decomposition theorem 
(see \cite[Section~4]{PU}).

The set $B_n(x,\e) := \{y \in X: d(f^ix,f^iy) < \e$ for $0 \le i < n\}$ is called a 
\textit{Bowen ball}.

A dynamical system $(X,f)$ is called \textit{conformal-like} if the image of every ball is a ball. 
A conformal map is a map
that  preserves angles and orientation; in the special case of smooth dynamics, the Jacobian of a conformal map is a positive multiple of a rotation matrice. Hence linear conformal maps preserve balls and are thus  conformal-like
but not every linear conformal-like map is conformal; for example it could reverse orientation.

\section{Proofs}
\subsection{Proof of Theorem~\ref{thm1}}
 From the locally eventually onto property, for each $y \in X$, and $\e> 0$ there is an $N(y,\e)\ge 1$ such that
$f^{N(y,\e)} (B(y,\e/3))= X$. Morover, by compactness of $X$ we can cover $X$ by a finite collection 
of balls $\{B(y_i,\e/3)\}_i$.  Let $N := \max_i\{N(y_i,\e)\}$.  Then since any ball $B(y,\e)$ contains one
of the $B(y_i,\e/3)$ we conclude that $f^N(B(y,\e))= X$ for all $y \in X$.

Now since $f$ is continuous and expanding, the image by $f^m$ of a Bowen ball $B_m(x,\e)$ is $B(f^m(x),\e)$, for every $0<\e<\delta_0$. 
Combining this with the previous paragraph yields $f^{m+N}B_m(x,\e) = X$ for each $x \in X$
and every $0<\e<\delta_0$.

Fix $\e > 0$ and choose $N$ as above. 
Consider an $N$-specification, i.e., a collection of orbit segments $\{f^{[a_j,b_j]} (x_j)\}_{j=1}^n$,
with $a_{i} - b_{i-1} \ge N$ for $2 \le i \le n$.
Setting $m_j := b_j - a_j$ and $N_j = m_j + N$
we have shown that 
$$
f^{N_j} (B_{m_j}( f^{a_j}(x_j),\e))  = X
\supset B_{m_{j+1}}( f^{a_{j+1}}(x_{j+1}),\e),$$
and thus
$$
B_{m_j}( f^{a_j}(x_j),\e)  
\cap f^{-N_j} ( B_{m_{j+1}}( f^{a_{j+1}}(x_{j+1}),\e) ) \neq \emptyset
$$
hold for each $1 \le j \le n$.
Iterating this, and noticing that $f$ is expanding,  yields that 
\begin{align}\label{chain}
\big\{ 
B_{m_1}( f^{a_1}(x_1),\e)  \cap  f^{- N_1 - N_2 - \cdots - N_i}  (B_{m_{i}}( f^{a_{i}}(x_{i}),\e)) 
\big\}_{2\le i \le n}
\end{align}
is a nested sequence of compact sets.
Any point in the intersection of these sets $\e$-shadows the specification,
and thus we have shown the specification property holds.

Finally we must show that the periodic specification property holds.
Fix $\e~>~0$ and consider an arbitrary $N$-specification
$\{f^{[a_j,b_j]} (x_j)\}_{j=1}^n$
 with $N$ choosen as above.
We extend this to a longer $N$ specification by choosing
$a_{n+1} = b_n + N$, $b_{n+1} = a_{n+1} + m_1$ and $x_{m+1} = f^{a_1 - a_{n+1}} x_1$.
Thus $m_{n+1} = m_1$ and 
$B_{m_{n+1}}(f^{a_{n+1}}x_{m+1},\e) = B_{m_{1}}(f^{a_{1}}x_{1},\e)$.
Therefore the  chain \eqref{chain} of containments extends to 
\begin{align*}
\big\{ 
B_{m_1}( f^{a_1}(x_1),\e)  \cap  f^{- N_1 - N_2 - \cdots - N_i}  (B_{m_{i}}( f^{a_{i}}(x_{i}),\e)) 
\big\}_{2\le i \le n+1}.
\end{align*}
The closure of the intersection of the extended chain of containments must contain a point 
fixed by $f^{N_1 + \cdots + N_{m+1}}$, hence the periodic specification property holds.
\hfill $\square$

\subsection{Proof of Theorem~\ref{prop:blokh}}

The strategy follows closely \cite[Appendix~A]{Bu}. For that reason we just give a brief sketch of the proof.
Let $X$ be a compact metric space and $f: X \to X$ be a continuous, locally eventually onto conformal map. 
The key step is a uniform control on the images of Bowen balls.
Indeed, while points in $n$-Bowen balls are within controlled distance to the original orbit during $n$ iterates, it is the size of the image 
the Bowen ball by iteration of $f^n$ which suggests how strong is the capability to obtain specification.  

\noindent {\bf Claim:} \textit{For any $\varepsilon>0$ there exists $\zeta(\varepsilon)>0$ so that} 
\begin{equation*}\label{big-image}
\mbox{diam}(f^n(B_n(x,\varepsilon))) \ge \zeta(\varepsilon)
\quad \text{\textit{for every}} \;n\ge 1\; \text{\textit{and}}\; x\in X.
\end{equation*}

\begin{proof}[Proof of the claim]
Fix $x\in X$. By conformality, for each $n\ge 1$ the set $f^{n}(B_{n}(x,\varepsilon))$ is a ball around $f^n(x)$.
Recall also that
\begin{equation}\label{e:n+1}
B_{n+1}(x,\varepsilon) = \bigcap_{j=0}^{n} f^{-j} (B(f^j(x),\e)) = B_{n}(x,\varepsilon) \cap f^{-n} (B(f^n(x),\e))
\end{equation}
and clearly
$
B_{n}(x,\varepsilon) \cap f^{-n} (B(f^n(x),\e)) \subseteq B_{n}(x,\varepsilon).
$

In particular,  by the conformality of $f$, for each $n\ge 1$ either: (i) the equality $B_{n+1}(x,\varepsilon) = B_{n}(x,\varepsilon)$ holds,
or (ii) the set $B_{n}(x,\varepsilon) \cap f^{-n} (B(f^n(x),\e))$ is strictly contained in $B_{n}(x,\varepsilon)$. 
In the second case, 
there exists a point $y\in B_{n}(x,\varepsilon)$ so that $f^n(y) \notin B(f^n(x),\e)$. 
This shows that 
the ball $f^{n}(B_{n}(x,\varepsilon)) \supset B(f^n(x),\varepsilon)$, combining with  \eqref{e:n+1} yields   
$$
f^{n} (B_{n+1}(x,\varepsilon)) = B(f^{n}(x),\varepsilon).
$$
Altogether, this proves that for every $n\ge 1$ there exists $0\le j < n $ so that 
$f^n(B_{n}(x,\varepsilon)) = f^{n-j} (B(f^j(x),\e))$.

Thus, in order to prove the claim it is enough to show that the forward image of balls of
a definite size do not degenerate: for any $\varepsilon>0$ there exists $\zeta(\varepsilon)>0$ such that 
$\mbox{diam}(f^n(B(z,\varepsilon))) \geq \zeta(\varepsilon)$ for every $n\ge 1$ and every $z\in X$.

Indeed, since $f$ is locally eventually onto, for any given $z\in X$ there exists $N(z,\varepsilon)>0$ such that 
$f^{N(z,\varepsilon)}(B(z,\varepsilon)) = X$; hence there exists 
$\zeta_z(\varepsilon)>0$ such that $\mbox{diam}(f^n(B(z,\varepsilon))) \geq \zeta_z(\varepsilon)$ for every $n\ge 1$.
The continuity of $f$ and compactness of $X$ ensures that $\min_{z\in X} \zeta_z(\varepsilon)>0$, proving the claim.
\end{proof}

We now claim that $f$ satisfies the periodic specification property.
Indeed, given $\e~>~0$ let $N=N(\varepsilon)\ge1$ be such that 
$f^N(B(x,\zeta(\varepsilon)))=X$ for every $x\in X$.
Such $N\ge 1$ does exists as $f$ is locally eventually onto and $X$ is compact.
The proof of the periodic specification property now follows as in Theorem~\ref{thm1}.
\hfill $\square$

\subsection{Proof of Theorem~\ref{thm-beta}}

Since items (ii) and (iii) are known (see e.g., \cite{Bu}) we need only prove that each $T_\beta$ is 
locally eventually onto.

Fix $\beta>1$ and take an arbitrary interval $J\subset [0,1)$. We claim that there exists $N\ge 1$ so that 
$T_\beta^N(J)=[0,1)$.
We may assume without loss of generality that $J$ is contained in some domain of smoothness for $T_\beta$.  By the mean value theorem, $\mbox{Leb}(T_\beta(J)) \ge \beta \,\mbox{Leb}(J)$. If  $T_\beta(J)\cap D_{T_\beta}=\emptyset$
then  $\mbox{Leb}(T_\beta^2(J)) \ge \beta^2 \, \mbox{Leb}(J)$. Since the diameter is bounded, a recursive argument shows
that $T_\beta^k(J)\cap D_{T_\beta}\neq \emptyset$ for some $k \ge 1$. In particular $T_\beta^k(J) \supset [0,a)$ for some $a\in (0,\frac1\beta]$.
Since $T_\beta(0)=0$, and $T_\beta$ is monotone increasing in  $[0,\frac1\beta]$ then there
exists $N\ge 1$  so that $T_\beta^N(J) \supset [0,\frac1\beta]$. This assures that  $T_\beta^{N+1}(J)=[0,1)$.
\hfill $\square$

\section{Examples}

We finish with some examples.
The first example is a simple examples of piecewise expanding continuous maps which need not be neither expansive nor transitive.

\begin{example}
Consider the continuous and piecewise expanding interval map $f_0: [0,\frac12] \to [0,\frac12]$
given by 
\[
f_0(x)
= 
\begin{cases}
\;\quad 3x & \; \mbox{if } x\in [0,\frac16] \\
-3x+1 & \; \mbox{if } x\in (\frac16,\frac13] \\
\;\; 3x-1 & \; \mbox{if } x\in (\frac13,\frac12].
\end{cases}
\]
Let $f:[0,1] \to [0,1]$ be obtained by replication of the dynamics $f_0$ in intervals of exponential decreasing growth accumulating $1$, 
defined by the relation 
$$
f(x)= 1-2^{-n} + 2^{-n} f_0(2^n (x-1+2^{-n})),
	\qquad x\in (1-2^{-n}, 1- 2^{-(n+1)}].
$$
and $f(0)=0$, $f(1)=1$.
Clearly $f$ is  piecewise expanding, continuous, not expanding nor transitive.
\end{example}

The next example shows that transitivity is essential to avoid unattainable repelling points.

\begin{example}
Consider the continuous and $C^1$-piecewise expanding interval map $f: [0,1] \to [0,1]$
given by 
\[
f(x)
= 
\begin{cases}
\;\quad 3x & \; \mbox{if } x\in [0,\frac13] \\
-2x+\frac53 & \; \mbox{if } x\in (\frac13,\frac23] \\
\;\;\; 2x-1 & \; \mbox{if } x\in (\frac23,1].
\end{cases}
\]
The map is not transitive as $f([\frac13,1]) = [\frac13,1]$, in other words, $[\frac13,1]$ is an $f$-invariant domain. 
Thus $f$ is not locally eventually onto. Nevertheless, the attractor $\Lambda:= \bigcap_{n\ge 0} f^n((0,1]])=[\frac13,1]$ 
and $f\mid_\Lambda$ is locally eventually onto.
\end{example}

Finally we complete this note with an example showing that locally eventually onto is weaker than specification. 
We consider an example suggested by Lindenstrauss (cf.\  \cite[Example~2.31]{AkAuNa}) of a locally eventualy onto
map having no periodic points.

\begin{example}\label{ex:Lindenstrauss}
Consider the subshift $Y_0\subset \{0,1,2\}^{\mathbb N}$ consisting of the set of sequences that admit no consecutive 0's, let
and let $\pi: Y_0 \to \{1,2\}^{\mathbb N}$ be given by supression of the 0's in the sequences belonging to $Y_0$. Endowing 
the shift spaces with the usual distances, $\pi$ is a continuous map on a compact metric space, hence it is uniformly continuous.

Consider a minimal subshift $X\subset  (\{1,2\}^{\mathbb N},\sigma)$ and let $Y=\pi^{-1}(X)$. Akin et al proved that $(Y,\sigma)$ 
is locally eventually onto (cf.\  Example~2.31 in \cite{AkAuNa}). We claim that $(Y,\sigma)$  does not satisfy the specification property.
Recall that a factor of a map 
 of a compact space with specification satisfies specification (cf.\  \cite[Proposition~21.4]{DGS}) 
This does not directly apply to our situation since we do not have compactness, however
 it is not hard to prove that the commuting diagram
$$
\begin{array}{ccc}
Y  & \rightarrow_{\sigma}& Y \\
\downarrow_\pi  & & \downarrow_\pi \\
X  & \rightarrow_{\sigma}& X \\
\end{array}
$$
together with the uniform continuity of $\pi$ ensures that if $(Y,\sigma)$ satisfies the 
specification property then so does $(X,\sigma)$. 
Second, $(X,\sigma)$ does not satisfy the specification property. Indeed, 
if $(X, T)$ has the specification property and its natural extension is expansive 
then $(X, T)$ has the periodic specification property (see e.g.,  Lemma 6 in \cite{KwLaOp}).
Altogether, this proves that $(Y,\sigma)$ does not satisfy the specification property, as claimed.
\end{example}

The following example, suggested by F. Przytycki, describes
an counter example to Theorem 1 if we do not assume that $X$ is a topological manifold. 

\begin{example}\label{Feliks}
Consider the circle $\mathbb S^1=\mathbb R/\mathbb Z$ and the doubling map $f: \mathbb S^1 \to \mathbb S^1$
given by $f(x)=2x \pmod{1}$. This is an expanding map (with constants $\lambda=2$ and $\delta_0=\frac12$), as
$d(f(x),f(y)) = 2 d(x,y)$ and $B(x,\frac12) \cap f^{-1}(z)$ is a singleton whenever $d(f(x),z) < \frac12$,
for all $x,y,z\in \mathbb S^1$.

Consider an enumeration $(p_n)_{n \in \N}$ of the set of periodic points for $f$, 
in such a way that their sequence of periods is non-decreasing, 
and choose a sequence
$(\eta_n)_{n \in \N}$ of positive real numbers converging quickly  to zero in such a way that
$$
K=\mathbb S^1\setminus \bigcup_{ n\ge 0} \bigcup_{k\ge 0} f^{-k}(B(p_n,\eta_n))
$$
is a non-empty compact subset of $\mathbb S^1$, and thus  $f(K) = K$. 

We furthermore suppose that the sequence $(\eta_n)_{n \in \N}$ is chosen as follows.
Since the periodic
points of $f$ are equidistributed in $\mathbb S^1$ (as $f$ is semi-conjugated to the full shift on two symbols),
for any $q_0 \in \mathbb S^1$ 
there exists $0<\zeta_0\ll 1$ so that 
$
\bigcup_{k\ge 0} f^{-k}(B(q_0,\zeta_0)) \cap \text{Per}(f) \neq \emptyset.
$
Set $q_0=p_0 \in \text{Per}(f)$.
As all points in $\mathbb S^1$ have dense pre-orbits we conclude that
$$K_0=\mathbb S^1\setminus \bigcup_{k\ge 0} f^{-k}(B(p_0,\zeta_0))
	 \supset K 
\neq \emptyset$$
is a Cantor set.  Let $n_1=\inf\{\ell \ge 1 \colon p_\ell \in K_0\}$ and write $q_1=p_{n_1}$. 
Choose $0<\zeta_1 \ll \zeta_0$
such that
$
\bigcup_{k\ge 0} f^{-k}\Big(B(q_0,\zeta_0)\cup B(q_1,\zeta_1)\Big) \cap \text{Per}(f) \neq \emptyset.
$
The previous condition can be assured by noting that any periodic point which intersects 
$B(q_0,\zeta_0)\cup B(q_1,\zeta_1)$ 
has combinatorics determined by either $q_0$ or $q_1$.
Then 
$$
K_1= \mathbb S^1\setminus \bigcup_{i= 0}^1 \bigcup_{k\ge 0} f^{-k}(B(q_i,\zeta_i))\subset K_0
$$
is a Cantor set which does not contain any of the periodic points in the set $\{p_n \colon 0\le n \le n_1\}$. 
Proceeding recursively,  we obtain a strictly decreasing sequence $(\zeta_\ell)_\ell$ of positive real numbers,
a strictly increasing sequence $(n_\ell)_\ell$ of positive integers and a nested sequence 
$(K_\ell)_\ell$ of Cantor sets such that 
$$
K_\ell =  \mathbb S^1\setminus \bigcup_{i= 0}^\ell \bigcup_{k\ge 0} f^{-k}(B(q_i,\zeta_i)) 
$$
contains some periodic point of $f$. Since the periodic points in $K_\ell\cap \text{Per}(f)$
are dense $K_\ell$ we have that $f(K_\ell)=K_\ell$ for every $\ell\ge 1$.
By construction, the set 
$$
K = \mathbb S^1\setminus \bigcup_{i \ge 0} \bigcup_{k\ge 0} f^{-k}(B(q_i,\zeta_i))  \neq \emptyset
$$
is a Cantor set having no periodic points and $f(K)=K$, as required.

Let us analyze the map $g :=f\mid_K$. This is clearly an expansive, and distance expanding map.
However, the following holds:
\begin{enumerate}
\item[(a)] $g$ has no periodic points;
\item[(b)] $g$ is not an open map; 
\item[(c)] $g$ is not an expanding map (i.e. condition (2) in Section~\ref{sec:def} fails).
\end{enumerate} 
Property (a) is immediate from the construction. Property (b) follows because for every 
open, distance expanding map, periodic points are dense in the non-wandering set 
(see Corollaries~4.2.4 and 4.2.5 in \cite{PU}). Property (c) is a consequence of property (b),
because of the equivalence between the notion of expanding  with  the notion of open, distance expanding 
on compact metric spaces (see \cite[Lemma~1]{CoRe}).
Moreover,
\begin{enumerate}
\item[(d)] $g$ is locally eventually onto.
\end{enumerate} 
This property is not immediate for subshifts (see e.g., Example~\ref{Petersen}).
In order to prove property (d) we will prove that each map $f\mid_{K_\ell}$ ($\ell\ge 1$) 
is locally eventually onto with uniform constants.

Fix any  $\ell \ge 1$ and $\varepsilon >0$. If $N=\lfloor \frac{-\log  \varepsilon}{\log 2}\rfloor$ then $f^N$
is a Markov map with $2^N$ full branches domains (injectivity domains), each of these with diameter larger than $2\varepsilon$. In particular, 
$f^N(B(x,\varepsilon))=\mathbb S^1$ for all $x\in \mathbb S^1$.
We claim that 
\begin{equation}\label{eq:approximation}
f^N\Big(B(x,\varepsilon) \cap K_\ell\Big) = K_\ell, \qquad \forall x \in \mathbb S^1.
\end{equation}
The inclusion $\subseteq$ is immediate. For the converse inclusion $\supseteq$, without loss of generality we can choose $\e = 2^{-N}$  and thus  $f^N\mid_{B(x,\varepsilon)}$ is a full branch for $f^N$ 
Since
$K_\ell$ has a dense set of the periodic points, every $z\in K_\ell$ is approximated by a sequence $(z_n)_n$ of 
periodic points in $K_\ell$. In particular, each of the points in the set  $f^{-N}(\{z\}) \cap B(x,\varepsilon)$ 
is an  accumulation point of periodic points in $K_\ell$.
Since $K_\ell$ is compact, this assures that $f^{-N}(\{z\})\cap B(x,\varepsilon) \in K_\ell$
and proves \eqref{eq:approximation}.
Now, as the Cantor sets are nested, one can use  ~\eqref{eq:approximation} to get 
$$
 f^N\Big(B(x,\varepsilon) \cap \bigcap_{\ell = 1}^n  K_\ell\Big) = 
\bigcap_{\ell =1}^n  K_\ell, \qquad \forall x \in \mathbb S^1, \; \forall n\ge 1.
$$
This, together with the continuity of the map $P(\mathbb S^1) \ni A \mapsto f^N(A)$ (in the Hausdorff topology)
implies that
$$
 f^N\Big(B(x,\varepsilon) \cap K\Big) = K, \qquad \forall x \in \mathbb S^1.
$$
This proves that $g=f\mid_K$ is locally eventually onto.

\end{example}

\begin{example}\label{Petersen}
There are examples of strongly mixing subshifts which are not locally eventually onto. Indeed, Petersen~ \cite{Pe} 
constructed a zero entropy, minimal and strongly mixing subshift $K\subset \{0,1\}^{\mathbb N}$
We claim that the distance expanding map $(K,\sigma)$ is not locally eventually onto. 

Assume, by contradiction that  $(K,\sigma)$ is locally eventually onto.  As $K$ is compact, for any $\varepsilon>0$ there exists $N=N(\varepsilon)\geq 1$
so that 
\begin{equation}
\label{eq:leo-balls}
\sigma^N(B_K(x,\varepsilon))=K \quad\text{ for every}\; x\in K.
\end{equation}
This implies that for any $0<\varepsilon < \frac12\text{diam} (K)$ there exist points $x_1,x_2 \in K$
so that $d(\sigma^N(x_1),\sigma^N(x_2))>\varepsilon.$ Hence, 
if $s(n,\varepsilon)$ denote the maximal cardinality of $(n,\varepsilon)$-separated subsets of $K$, a recursive
argument using \eqref{eq:leo-balls} 
together with the observation that $\sigma^n(B(x,n,\varepsilon))=B(\sigma^n(x),\varepsilon)$ for every $x\in K$, $n\ge 1$ and $\varepsilon>0$
ensures that $s(kN,\varepsilon)\geq 2^k$ for every $k\ge 1$.  Hence
\begin{align*}
h_{\text{top}}(\sigma\mid_K) 
	&=\lim_{\varepsilon \to 0} \limsup_{n\to\infty} \frac1n \log s(n,\varepsilon) \\
	& \geq \lim_{\varepsilon \to 0} \limsup_{k\to\infty} \frac1{kN} \log s(kN,\varepsilon)  \geq \frac1N \log 2 >0,
\end{align*}
which leads to a contradiction. This proves that $(K,\sigma)$ is not locally eventually onto. 
\end{example}

The next simple example illustrates that compactness is an essential assumption in Theorem~\ref{thm1}.

\begin{example}\label{non-compact}
The mixing and specification properties have been extensively studied in the case of symbolic dynamics (see e.g., \cite[Section~8]{KwLaOp}
and references therein). Here we give an example of a shift space, hence distance expanding, which is locally eventually
onto, has dense periodic orbits but for which the specification property fails. 

Consider the subshift 
$\Sigma_{\mathcal G}\subset \mathbb N^{\mathbb N}$ determined by the countable graph $\mathcal G$
with countable states $\mathbb N$ and whose allowed directed paths $v\to w, \; v,w\in \mathbb N$ are 
$0 \to w$ for every $w\in \mathbb N$, and the arrows $v \to w$ with $v\neq 0$ are admissible if and only if 
$w\in \{v-1,v\}$. The cylinder sets are defined by 
$
[v_0, v_1, \dots, v_n]=\big\{(w_0, w_1, w_2, \dots) \colon w_i=v_i, \, \forall 0\le i \le n\big\}.
$ 
The shift $\sigma: \Sigma_{\mathcal G} \to \Sigma_{\mathcal G}$ is locally eventually onto 
because 
$$
\sigma^{n+1}([v_0, v_1, \dots, v_n])\supset \sigma([1])=\Sigma_{\mathcal G}\quad 
\text{for every cylinder} \; [v_0, v_1, \dots, v_n].
$$
It is a simple exercise to show that the condition $\sigma^j([n])\cap [1]=\emptyset$ for every $0\le j \le n-1$
is incompatible with the specification property.
\end{example}

In final example we 
prove an optimality of Theorem~\ref{thm1}, in the sense that it fails if condition (2) in the definition of expanding map is removed. 

\begin{example}\label{embedding}
Let $\sigma: \{0,1\}^{\mathbb N} \to \{0,1\}^{\mathbb N}$ be the full shift, and let $g: [0] \to [0]$ be the first return map
of $\sigma$ to the cylinder $[0]$. More precisely, if $\tau : [0] \to \mathbb N$ is the first return time to $[0]$ given by
$$
\tau(x_0,x_1,x_2, x_3, \dots) = \inf\{ k\ge 1 \colon x_k=0 \}
$$
then $g(\cdot)=\sigma^{\tau(\cdot)}(\cdot)$. Equivalently, if $x_0=0$ then 
$$
g(x_0,x_1, x_2, x_3, x_4, \dots) = (x_k, x_{k+1}, x_{k+2}, \dots) 
$$ 
where $k=\tau(x_0,x_1,x_2, x_3, \dots)$.
After identification of the set $\mathcal S$ of all finite words 
 $(0, 1, 1, 1, \dots, 1,  0)$ with the cylinder $[0, 1, 1, 1, \dots, 1,  0] \subset \{0,1\}^{\mathbb N}$,
the map $g$ acts as a full shift $\mathcal S^{\mathbb N}$. By the previous identification, we 
will consider $\mathcal S^{\mathbb N}$ as a subset of the cylinder $[0]\subset \{0,1\}^{\mathbb N}$. This construction is often called the ``Rome graph".

Let $\Sigma \subset \mathcal S^{\mathbb N}$ be the locally eventually onto subshift so that
$g\mid_{\Sigma}:\Sigma \to\Sigma$ does not satisfy the specification property induced by  Example~\ref{non-compact}. 
Indeed, just use the bijection 
\begin{equation}\label{eq:bijection}
\mathbb N \to \mathcal S \quad\text{given by}\quad n \mapsto (0, \underbrace{1,1, \dots, 1}_{n}, 0)
\end{equation}
to embed the subshift $\Sigma_{\mathcal G}\subset \mathbb N^{\mathbb N}$ onto such subshift $\Sigma\subset \mathcal S^{\mathbb N}$. The possible unbounded amount of 1's in ~\eqref{eq:bijection} makes the subshift $\Sigma\subset \mathcal S^{\mathbb N}$ not closed.  Now, consider the $\sigma$-invariant and compact set 
$K\subset \{0,1\}^{\mathbb N}$ obtained as the closure of the saturated set
$$
{\bigcup_{n\ge 1} \bigcup_{j=0}^{n-1} \sigma^j\big(\big\{w\in \Sigma \colon \tau(w)=n \big\}\big)}.
$$
By construction and the fact that $K$ is closed we get  
$$
K\cap [0] = \text{closure}(\mathcal S^{\mathbb N}) = \Sigma \cup \{01^\infty\}
$$
and
\begin{align*}
K\cap [1] 
& = \text{closure}\Big({\bigcup_{n\ge 1} \bigcup_{j=1}^{n-1} \sigma^j\big(\big\{w\in \Sigma \colon \tau(w)=n \big\}\big)}\Big) \\
& = \Big({\bigcup_{n\ge 1} \bigcup_{j=1}^{n-1} \sigma^j\big(\big\{w\in \Sigma \colon \tau(w)=n \big\}\big)}\Big) 
	\cup \{1^\infty\}
\end{align*}

Note that the elements in $\{01^\infty, 1^\infty\}$ do not return to the cylinder $[0]$. Thus, using that $K\setminus \{01^\infty, 1^\infty\}$  is obtained by the union of the finite pieces of orbits of points in $\Sigma$ until their
first return time to $\Sigma$,  the distance expanding map $(K,\sigma\mid_K)$ does not satisfy the specification property.

We claim that $(K,\sigma)$ is locally eventually onto. Dealing with the induced topology, it is enough to prove that 
for any cylinder 
$[x_1,x_2, \dots, x_n]\cap K\neq\emptyset$ there exists $N\ge 1$ so that $\sigma^N([x_1,x_2, \dots, x_n])=\{0,1\}^{\mathbb N}$.
This is a consequence of the fact that $g$ is a Poincar\'e first return map of $\sigma$ to the global cross-section 
$[0]$ and that $g$ is locally eventually onto (recall Example~\ref{non-compact}).
Then Theorem~\ref{thm1} implies that the long inverse branches condition (2) in the definition of expanding map fails.
\end{example}


\begin{thebibliography}{99}
\bibitem{AkAuNa} Akin, Ethan; Auslander, Joseph; Nagar, Anima.
\textit{Variations on the concept of topological transitivity},
Studia Math.\ 235 (2016), no. 3, 225--249.

\bibitem{BM}
Barge, Marcy; Martin, Joe.
\textit{Dense orbits on the interval.}
Michigan Math.\ J.\ 34 (1987), no. 1, 3--11. 

\bibitem{Bl}
Blanchard, Fran\c cois. 
\textit{$\beta$-expansions and symbolic dynamics. }
Theoretical Computer Sci.\ 65 (1989), no. 2, 131--141.

\bibitem{Blo0}
Blokh, Alexander.
\textit{On Transitive Maps of One-Dimensional Branched Manifolds}.
Differential - Difference Equations and Problems of Mathematical Physics, Kiev (1984), 3--9.

\bibitem{Blo}
Blokh, Alexander.
\textit{The Spectral Decomposition for One-Dimensional Maps}
Dynamics Reported 4 (1995), 1--59, Springer, Berlin, 1995.  
 
 \bibitem{Bo}
 Bowen, Rufus.  \textit{Periodic points and measures for Axiom A diffeomorphisms}, Trans.\ Amer.\ Math.\ Soc.\ 154 (1971), 377--397.

\bibitem{Bu} Buzzi, J\'er\^ome.
\textit{Specification of the interval},
Trans. Amer. Math.\ Soc.\ 349 (1997) 2737--2754.

\bibitem{CM}
 Coven, Ethan; Mulvey,  Irene. 
 \textit{Transitivity and the centre for maps of the circle.} Ergod.\ Th.\ Dynam.\ Sys.\ 6 (1986), no. 1, 1--8.

\bibitem{CoRe} Coven, Ethan; Reddy, William.
\textit{Positively expansive maps of compact manfolds}
 In: Nitecki Z., Robinson C. (eds) Global Theory of Dynamical Systems. Lecture Notes in Mathematics, vol 819. Springer, Berlin, Heidelberg

\bibitem{DGS} 
Denker, Manfred; Grillenberger, Christian; Sigmund, Karl.
 Ergodic theory on compact spaces. Lecture Notes in Mathematics, Vol. 527. Springer-Verlag, Berlin-New York, 1976. iv+360 pp.


\bibitem{KH} Katok, Anatole; Hasselblatt, Boris.
\newblock \textit{Introduction to the Modern Theory of Dynamical Systems},
\newblock Cambridge University Press, 1995.

\bibitem{KKO} Kulczycki, Marcin; Kwietniak, Dominik; Oprocha, Piotr, 
\newblock \textit{On almost specification and average shadowing properties},
Fund.\ Math.\ 224 (2014), no. 3, 241--278.

\bibitem{KwLaOp} Kwietniak, Dominik, {\L}acka, Martha, Oprocha, Piotr.
\textit{A panorama of specification-like properties and their consequences} 
 Dynamics and Numbers, Contemporary Mathematics, vol. 669, (2016) 155--186.

\bibitem{Ol}
Oliveira, Krerley. 
\textit{Every expanding measure has the nonuniform specification property.}
Proc.\  Amer.\  Math.\ Soc.\ 140:4 (2006) 1309--1320.

\bibitem{Pe} 
Petersen, Karl. 
\textit{A Topologically Strongly Mixing Symbolic Minimal Set.}
Trans.\ Amer.\ Math.\ Soc.\ 148:2  (1970) 603--612.

\bibitem{PU} 
Przytycki, Feliks; Urba\'nki Mariusz.
\textit{Conformal Fractals: Ergodic Theory Methods,} 
 London Mathematical Society Lecture Note Series. 371,
 Cambridge University Press, Cambridge, 2010. x+354 pp.
 
\bibitem{Ru}  Ruette, Sylvie. 
\textit{Chaos on the interval}. 
University Lecture Series, 67. American Mathematical Society, Providence, RI, 2017. xii+215 pp

\bibitem{Sc} 
 Schmeling, J\"org.
\textit{Symbolic dynamics for $\beta$-shifts and self-normal numbers}.
Ergod.\ Th.\ Dynam.\ Sys.\ 17:3 (1997)  675--694.
 
 \bibitem{Va}  
Varandas, Paulo. 
\textit{Non-uniform specification and large deviations for weak Gibbs measures,} 
J.\ Stat.\ Phys., 146 (2012) 330--358.
 
\bibitem{YaYiWa} Yan, Qi; Yin, Jiandong; Wang, Tao.
\textit{Some weak specification properties and strongly mixing.}
Chinese Annals Mathematics\ Ser.\ B 38 (2017), no. 5, 1111--1118.

\end{thebibliography}
\end{document}